# Transcendental Numbers and the Lambert-Tsallis Function

J. L. E. da Silva          R. V. Ramos

leonardojade@alu.ufc.br   rubens.ramos@ufc.br

*Lab. of Quantum Information Technology, Department of Teleinformatic Engineering – Federal University of Ceara - DETI/UFC, C.P. 6007 – Campus do Pici - 60455-970 Fortaleza-Ce, Brazil.*

*Abstract*

To decide upon the arithmetic nature of some numbers may be a non-trivial problem. Some cases are well know, for example exp(1) and $W$(1), where $W$ is the Lambert function, are transcendental numbers. The Tsallis $q$-exponential, $e_q(z)$, and the Lambert-Tsallis $W_q(z)$ function, where $q$ is a real parameter, are, respectively, generalizations of the exponential and Lambert functions. In the present work we use the Gelfond-Schneider theorem in order to show the arithmetic conditions on $q$ and $z$ such that $W_q(z)$ and $\exp_q(z)$ are transcendental.

## 1. Introduction

To determine if a given number is algebraic irrational or transcendental may not be an easy task. For example, up to the moment there is no answer for the possible transcendence of $e + \pi$, $e^{\pi}$ and $\log(\pi)$. On the other hand, it is well know that $e$ and $W$(1), where $W$ is the Lambert function, are transcendental number. This immediately gives rise to the question if the functions $e^z$ and $W(z)$ are good producers of transcendental numbers. The answer is positive for the exponential but, up to now, it is not known the answer for $W(z)$. On the other hand, one can check if their deformed versions, the Tsallis $q$-exponential, $e_q(z)$, and the Lambert-Tsallis $W_q(z)$ function are good producers of transcendental numbers. In the present work, we use the Gelfond-Schneider theorem to answer positively this question.

This work is outlined as follows: In Section 2, the Tsallis $q$-exponential and the Lambert-Tsallis function are briefly reviewed. In Section 3 we use the Gelfond-Schneider theorem to determine the arithmetic conditions that should be obeyed by $q$ and $z$ in order to $exp_q(z)$ and $W_q(z)$ to be transcendental numbers. At last, the conclusions are drawn in Section 4.

## 2. Tsallis $q$-exponential and Lambert-Tsallis $W_q$ function

The Lambert $W$ function is a well-known elementary function that has been used in different areas of mathematics, computer Science and physics [1-6]. Basically, the Lambert $W$ function is the solution of the equation

$$W(z)e^{W(z)} = z. \tag{1}$$

On the other hand, if one uses the $q$-exponential function proposed by Tsallis [7]

$$e_q^z = \begin{cases} e^z & q=1 \\ \left[1+(1-q)z\right]^{1/(1-q)} & q \neq 1 \ \& \ 1+(1-q)z \geq 0 \\ 0 & q \neq 1 \ \& \ 1+(1-q)z < 0 \end{cases} \tag{2}$$

in the Lambert equation (1), one obtains the Lambert-Tsallis equation

$$W_q(z)e_q^{W_q(z)} = z, \tag{3}$$

whose solutions are the Lambert-Tsallis $W_q$ functions [8]. More details about the $W_q$ function can be found in [8-13]. It can be found that the branch point of the Lambert-Tsallis $W_q$ function is ($z_b = \exp_q(1/(q-2))/(q-2)$, $W_q(z_b) = 1/(q-2)$), for $q \neq 2$. There is no branch point with finite $z_b$ for $q = 2$. For $q = 1$, the branch point of Lambert $W$ function $(-1/e,-1)$ is recovered. The real solution in the interval $z_b \leq z < 0$ is $W_q^-(z)$ while the real solution in the interval $z_b \leq z < \infty$ is $W_q^+(z)$. The function $W_q^+(z)$ keeps its concavity according to $d^2W_q^+(z)/dz^2 < 0$. On the other hand, $W_q^-(z)$ decreases from the branch point and goes toward the point $(0^-, -\infty)$.

The first derivatives of $W_q$ is given by

$$\frac{dW_q(z)}{dz} = \frac{-\left((1-q)W_q(z)+1\right)^{\frac{q}{q-1}}}{(q-2)W_q(z)-1}. \tag{4}$$

The coordinates of the branches points of the functions $W_q(z)$ are obtained doing $dW_q/dz = \infty$ (points of the curve with vertical tangent). The first solution is $W_q = -\infty$. The second solution, $W_q = W_q^b$, depends on the value of $q$, as shown before.

## 3. Transcendentality of $W_q(x)$ and $\exp_q(x)$

Firstly, we consider the number $W_q(1)$. Using eq. (2) in eq. (3) and $z = 1$ one gets the following polynomial

$$x^{q-1} - (1-q)x - 1 = 0 \Rightarrow x^{q-1} = (1-q)x + 1, \tag{5}$$

where $x = W_q(1)$. The Gelfond-Schneider theorem states that, if $\alpha$ and $\beta$ are algebraic numbers, with $\alpha \neq 0$ or $1$, and $\beta$ irrational, then $\alpha^\beta$ is transcendental [14].

**Theorem 1** – If $q$ is an algebraic and irrational number, then $W_q(1)$ is transcendental.

Proof: Initially $q$ is an algebraic and irrational number. Let us assume that $x$ in (5) is algebraic. In this case, according to Gelfond-Schneider theorem $x^{q-1}$ is transcendental. On the other hand, since $q$ and $x$ are algebraic, $(1-q)x+1$ is also algebraic, since the set of algebraic numbers is a field. Hence, there is a contradiction since the left side of (5) is transcendental and the right side of eq. (5) is algebraic. Thus, $x = W_q(1)$ must be transcendental □.

The transcendentality of $e_q(z)$ is explained in Theorem 2:

**Theorem 2** – If $q$ is an algebraic and irrational number and $z$ ($\neq 0$) is algebraic then $e_q(z)$ is transcendental.

Proof: Using eq. (2) one has

$$e_q^z = \left[1+(1-q)z\right]^{\frac{1}{(1-q)}}. \tag{6}$$

Since $q$ and $z$ are algebraic numbers, $1+(1-q)z$ is also algebraic. Similarly, since $q$ is algebraic and irrational, $1/(1-q)$ is also algebraic and irrational. Thus, according to Gelfond-Schneider theorem $e_q(z)$ is transcendental □.

At last, the transcendentality of $W_q(z)$ is discussed in Theorem 3:

**Theorem 3** – If $q$ is an algebraic and irrational number and $z$ ($\neq 0$) is algebraic then $W_q(z)$ is transcendental.

Proof: Using eq. (3) one has

$$W_q(z) e_q^{W_q(z)} = z \Rightarrow e_q^{W_q(z)} = \frac{z}{W_q(z)} \tag{7}$$

Let us assume that $q$ is algebraic and irrational and $W_q(z)$ is algebraic. In this case, according to Theorem 2 $\exp_q(W_q(z))$ is transcendental. On the other hand, $z/W_q(z)$ is the ration between two algebraic numbers, that is also an algebraic number. Therefore, there is a contradiction and $W_q(z)$ cannot be algebraic □.

## 4. Other results using $W_q(z)$, $\exp_q(z)$ and $\log_q(z)$

The first derivative of functions that generate transcendental numbers for algebraic arguments may not always generate transcendental image. For example, the derivative of $e^z$ is itself, and it is known that for non-null algebraic arguments, the image of this function is always a transcendental number. However, the same does not happen with the function *log(z)* whose first derivative is 1/*z*, which for non-null algebraic arguments always assumes algebraic values. In this perspective, we present the following theorem

**Theorem 4** – If $z_0$ ($\neq 0$) is algebraic and $q$ is an algebraic and irrational number, then $d(log_q(z_0))/dz$ is transcendental, where $log_q(z)$ is the Tsallis *q*-logarithmic function.

Proof: Using the derivative of the $ln_q(z_0)$ function

$$\frac{d}{dz}\left(\ln_q(z)\right) = \frac{1}{z^q} \quad (z > 0, \forall q), \tag{8}$$

Thus, for $z = z_0$ non-null algebraic and $q$ algebraic and irrational, by the Gelfond-Schneider theorem one has that $1/(z_0)^q$ is transcendental ☐.

It is known that $e^\pi$ is a transcendental number, since $e^\pi = (-1)^{-i}$. However, nothing is known about the transcendentality of $e^e$. For the *q*-exponential function one has the following result

**Theorem 5** – If $z$ is transcendental and $q \in \mathbb{Q}\backslash\{1\}$. Then $\exp_q(z)$ is transcendental.

Proof: Since a transcendental number power to a rational number is transcendental, from eq. (6), $\exp_q(z)$ with $z$ transcendental and $q \in \mathbb{Q}\backslash\{1\}$ is also transcendental ☐.

Theorem 5 allows us to show that $e_q^e$ and $e_q^\pi$, for $q \in \mathbb{Q}\backslash\{1\}$, are also transcendental numbers.

Some values of *q* the Lambert-Tsallis $W_q$ function return algebraic values for non-null algebraic argument, for example, $W_2(1) = ½$ [8]. So, for these cases we introduce the following theorem

**Theorem 6** – If $z$ is non-null algebraic and $W_q(z) \in \mathbb{Q}\setminus\{\mathbb{Z}\}$. Then $W_q(z)^{W_q(z)^{W_q(z)}}$ is transcendental.

Proof: If $W_q(z) \in \mathbb{Q}\setminus\{\mathbb{Z}\}$, then $W_q(z)^{W_q(z)}$ is irrational algebraic [15]. So at, by Gelfond-Schneider theorem, $W_q(z)^{W_q(z)^{W_q(z)}}$ is transcendental □.

## 4. Conclusions

If $q$ is an algebraic and irrational number and $z$ is a not null algebraic number, then $W_q(z)$ and $\exp_q(z)$ are transcendental numbers. Hence, the $q$-exponential function and the Lambert-Tsallis function are good producers of transcendental numbers.

## Acknowledgements

This study was financed in part by the Coordenação de Aperfeiçoamento de Pessoal de Nível Superior - Brasil (CAPES) - Finance Code 001, and CNPq via Grant no. 307184/2018-8. Also, this work was performed as part of the Brazilian National Institute of Science and Technology for Quantum Information.